\documentclass[11pt]{article}
\usepackage{graphicx}
\usepackage{amsfonts}
\usepackage{latexsym}
\usepackage{epsfig}

\newcommand{\ignore}[1]{}
\newcommand{\beql}[1]{\begin{equation}\label{#1}}
\newcommand{\eeq}{\end{equation}}

\newcommand{\eqref}[1]{{\rm (\ref{#1})}}

\newcommand{\Abs}[1]{{\left|{#1}\right|}}

\newcommand{\Qed}{\ \\\mbox{$\Box$}}

\newcommand{\Set}[1]{{\left\{{#1}\right\}}}

\newcommand{\RR}{{\mathbb R}}

\newcommand{\ZZ}{{\mathbb Z}}

\newcommand{\dens}{{\rm dens\,}}

\newcommand{\supp}{{\rm supp\,}}

\newcommand{\ft}[1]{\widehat{#1}}

\def\eps{\epsilon}

\def\rr{{\RR}}
\def\zz{{\ZZ}}

\def\half{\frac{1}{2}}

\newcounter{open}
\newcounter{dfn}
\def\thedfn{\arabic{dfn}}

\newcounter{obs}
\def\theobs{\arabic{obs}}

\newcounter{thm}
\newcounter{mysec}
\newcounter{mysubsec}[mysec]
\newtheorem{theorem}{Theorem}
\newtheorem{corollary}{Corollary}
\newtheorem{lemma}{Lemma}

\begin{document}

\title{Tiling and spectral properties of near-cubic domains}
\author{Mihail N. Kolountzakis and Izabella {\L}aba}
\maketitle


\section{Introduction}


Let $E$ be a measurable set in $\RR^n$ such that $0<|E|<\infty$.
We will say that $E$ {\em tiles $\RR^n$ by translations} if there is a discrete
set $T\subset\RR^n$ such that, up to sets of measure 0, the sets $E+t:\
t\in T$ are mutually disjoint and $\bigcup_{t\in T}(E+t)=\RR^n$. We
call any such $T$ a {\em translation set} for $E$, and write 
$E+T=\RR^n$.
A tiling $E+T=\RR^n$ is called {\em periodic} if it admits a 
period lattice of rank $n$; it is a {\em lattice tiling} if
$T$ itself is a lattice.  Here and below, a {\em lattice}
in $\RR^n$ will always be a set of the form $T\ZZ^n$, where
$T$ is a linear transformation of rank $n$.

It is known \cite{V}, \cite{M} that if a convex set $E$ tiles $\RR^n$ by
translations, it also admits a 
lattice tiling.  A natural question is whether a similar result holds if $E$
is ``sufficiently close" to being convex, e.g. if it is close enough (in
an appropriate sense) to a $n$-dimensional cube.  In this paper we prove
that this is indeed so in dimensions 1 and 2; we also construct a
counterexample in dimensions $n\geq 3$.

A major unresolved problem in the mathematical theory of tilings is
the {\em periodic tiling conjecture}, which
asserts that any $E$ which tiles $\rr^n$ by translations must also admit
a periodic tiling.  (See \cite{GS} for an overview of this and other related
questions.) The conjecture has been proved for all bounded measurable
subsets of $\RR$ \cite{LW1}, \cite{KL} and
for topological discs in $\RR^2$ \cite{GN}, \cite{Ke}.
Our Theorem \ref{thm.2d} and Corollary \ref{cor.2d} prove
the conjecture for near-square domains in $\RR^2$.  We emphasize that
no assumptions on the topology of $E$ are needed; in particular, 
$E$ is not required to be connected and may have infinitely
many connected components.

Our work was also motivated in part by a conjecture of Fuglede
\cite{Fug}.  We call a set $E$ {\em spectral} if there is a
discrete set  $\Lambda\subset\RR^n$, which we call a {\em spectrum} for
$E$, such that $\{e^{2\pi i \lambda\cdot x}:
\lambda\in \Lambda\}$ is an orthogonal basis for $L^2(E)$.  Fuglede
conjectured that $E$ is spectral if and only if it tiles $\RR^n$ by
translations, and proved it under the assumption that either the
translation set $T$ or the spectrum $\Lambda$ is a lattice.  This
problem was addressed in many recent papers (see e.g. \cite{IKP}, \cite{JP1},
\cite{K1}, \cite{KP}, \cite{L1}, \cite{L2}, \cite{LW1}, \cite{LW2}), and in
particular the conjecture has been proved for convex regions in $\rr^2$
\cite{K3}, \cite{IKT1}, \cite{IKT2}.

It follows from our Theorem \ref{th:main-1dim} and from
Fuglede's theorem that the conjecture is true for $E\subset\RR$ such that
$E$ is contained in an interval of length strictly less than $3|E|/2$.
(This was proved in \cite{L2}
in the special case when $E$ is a union of finitely many intervals
of equal length.)  In dimension 2, we obtain the ``tiling $\Rightarrow$
spectrum" part of the conjecture for near-square domains.  Namely,
if $E\subset\rr^2$ tiles $\RR^2$ and satisfies the assumptions of Theorem
\ref{thm.2d} or Corollary \ref{cor.2d}, it also admits a lattice tiling,
hence it is a spectral set by Fuglede's theorem
on the lattice case of his conjecture. We do not know how to
prove the converse implication.

Our main results are the following.

\begin{theorem}\label{th:main-1dim}
Suppose $E\subseteq[0,L]$ is measurable with measure $1$ and
$L = 3/2-\epsilon$ for some $\epsilon>0$.
Let $\Lambda\subset\RR$ be a discrete set containing $0$.
Then\\
(a) if $E+\Lambda =\RR$ is a tiling, it follows that $\Lambda=\ZZ$.\\
(b) if $\Lambda$ is a spectrum of $E$, it follows that $\Lambda=\ZZ$.
\end{theorem}

The upper bound $L<3/2$ in Theorem \ref{th:main-1dim} is optimal:
the set $[0,1/2]\cup [1,3/2]$ is contained in an interval of length
$3/2$, tiles $\ZZ$ with the translation set $\{0,1/2\}+2\zz$, and has
the spectrum $\{0,1/2\}+2\zz$, but does not have either a lattice
translation set or a lattice spectrum.  This example has been
known to many authors; an explicit calculation of the spectrum is
given e.g. in \cite{L1}.

\begin{theorem}\label{thm.2d}
Let $E\subset \RR^2$ be a measurable set such that
$[0,1]^2\subset E\subset [-\eps,1+\eps]^2$
for $\eps>0$ small enough.
Assume that $E$ tiles $\RR^2$ by translations.  Then $E$ also admits
a tiling with a lattice $\Lambda\subset\RR^2$ as the translation set.
\end{theorem}

Our proof works for $\eps<1/33$; we do not know what is the optimal
upper bound for $\eps$.

\begin{figure}[htbp] \centering
\ \psfig{figure=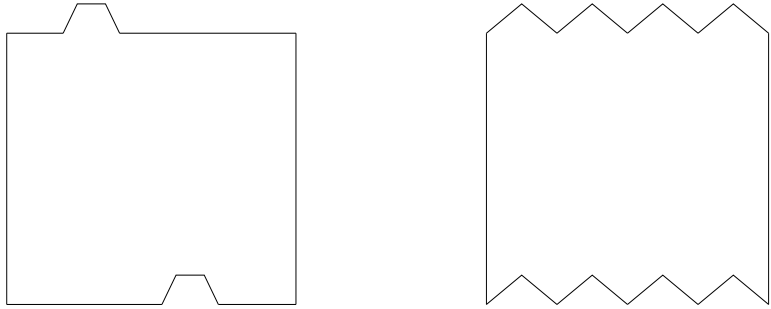,height=1.6in,width=4in}
\caption{Examples of near-square regions which tile $\RR^2$.  Note that
the second region also admits aperiodic (hence non-lattice) tilings.}
        \label{fig:squares}
        \end{figure}

\begin{corollary}\label{cor.2d}
Let $E\subset \RR^2$ be a measurable set such that
$|E|=1$ and $E$ is contained in a square of sidelength $1+\eps$ for
$\eps>0$ small enough. If $E$ tiles $\RR^2$ by translations, then it
also admits a lattice tiling.
\end{corollary}

\begin{theorem}\label{3d-thm}
Let $n\geq 3$. Then for any $\eps>0$ there is a set $E\subset\RR^n$ 
with $[0,1]^n\subset E\subset [-\eps,1+\eps]^n$ such that $E$ tiles
$\RR^n$ by translations, but does not admit a lattice tiling.
\end{theorem}


\section{The one-dimensional case}


In this section we prove Theorem \ref{th:main-1dim}.
We shall need the following crucial lemma.

\begin{lemma}\label{lm:shifts}
Suppose that $E\subseteq[0,L]$ is measurable with measure $1$ and
that $L = 3/2-\epsilon$ for some $\epsilon>0$.
Then
\beql{overlaps}
\Abs{E \cap (E+x)} > 0\ \ \mbox{whenever $0\le x < 1$}.
\eeq
\end{lemma}

{\bf Proof of Lemma \ref{lm:shifts}.}
We distinguish the cases (i) $0\le x \le 1/2$, (ii) $1/2 < x \le 3/4$ and
(iii) $3/4 < x < 1$.

\noindent
\underline{(i) $0\le x \le 1/2$}

This is the easy case as $E \cup (E+x) \subseteq [0, L+1/2] = [0, 2-\epsilon]$.
Since this interval has length less than $2$, the sets $E$ and $E+x$ must intersect in
positive measure.

\noindent
\underline{(ii) $1/2 < x \le 3/4$}

Let $x = 1/2 + \alpha$, $0< \alpha \le 1/4$.
Suppose that $\Abs{E \cap (E+x)} = 0.$
Then $1+2\alpha \le 3/2$ and
$$
\Abs{(E \cap [0, x]) \cup (E \cap [x, 2x])} \le x,
$$
as the second set does not intersect the first when shifted back by $x$.
This implies that
$$
\Abs{E} \le x + (3/2 - \epsilon - 2x) = 3/2 - \epsilon - x = 1 - \epsilon - \alpha < 1,
$$
a contradiction as $\Abs{E} = 1$.

\noindent
\underline{(iii) $3/4 \le x < 1$}

Let $x = 3/4 + \alpha$, $0<\alpha < 1/4$.
Suppose that $\Abs{E \cap (E+x)} = 0.$
Then
$$
\Abs{(E \cap [0, 3/4-\alpha-\epsilon]) \cup (E \cap [3/4 + \alpha, 3/2-\epsilon])}
 \le 3/4 - \alpha - \epsilon,
$$
for the second set translated to the left by $x$ does not intersect the first.
This implies that
$$
\Abs{E} \le (3/4 - \alpha - \epsilon) + 2\alpha + \epsilon = 3/4 + \alpha < 1,
$$
a contradiction.
\Qed

We need to introduce some terminology.
If $f$ is a nonnegative integrable function on $\RR^d$ and $\Lambda$
is a subset of $\RR^d$, we say that $f+\Lambda$ is a packing
if, almost everywhere,
\beql{packing}
\sum_{\lambda\in\Lambda} f(x-\lambda) \le 1.
\eeq
We say that $f+\Lambda$ is a tiling if equality holds almost everywhere.
When $f = \chi_E$ is the indicator function of a measurable
set, this definition coincides with the classical
geometric notions of packing and tiling.

We shall need the following theorem from \cite{K1}.

\begin{theorem}\label{th:tiling-from-packing}
If $f, g \ge 0$, $\int f(x) dx = \int g(x) dx = 1$ and
both $f+\Lambda$ and $g+\Lambda$ are packings of $\RR^d$, then
$f+\Lambda$ is a tiling if and only if $g+\Lambda$ is a tiling.
\end{theorem}

\noindent
{\bf Proof of Theorem \ref{th:main-1dim}.}
\noindent (a)
Suppose $E + \Lambda$ is a tiling.  From Lemma \ref{lm:shifts} it follows
that any two elements of $\Lambda$ differ by at least $1$.
This implies that $\chi_{[0,1]}+\Lambda$ is a packing,
hence it is also a tiling by Theorem \ref{th:tiling-from-packing}.
Since $0\in\Lambda$, we have $\Lambda = \ZZ$.

\noindent (b)
Suppose that $\Lambda$ is a spectrum of $E$.
Write
$$
\delta_\Lambda = \sum_{\lambda\in\Lambda}\delta_\lambda
$$
for the measure of one unit mass at each point of $\Lambda$.
Our assumption that $\Lambda$ is a spectrum for $E$ implies that
$$
\Abs{\ft{\chi_E}}^2 + \Lambda = \RR
$$
is a tiling (see, for example, \cite{K1}).
This, in turn, implies that $\dens \Lambda = 1$.

We now use the following result from \cite{K1}:
\begin{theorem}\label{th:tiling-implies-supp}
Suppose that $f\ge 0$ is not identically $0$, that $f \in L^1(\RR^d)$,
$\widehat{f}\ge 0$ has compact support and $\Lambda\subset\RR^d$.
If $f+\Lambda$ is a tiling then
\beql{supp-cond-1}
\supp \ft{\delta_\Lambda} \subseteq \Set{\ft{f} = 0} \cup \Set{0}.
\eeq
\end{theorem}

Let us emphasize here that the object $\ft{\delta_\Lambda}$,
the Fourier Transform of the tempered measure $\delta_\Lambda$,
is in general a tempered distribution and need not be a measure.

For $f = \Abs{\ft{\chi_E}}^2$ Theorem \ref{th:tiling-implies-supp} implies
\beql{support}
\supp{\ft{\delta_\Lambda}} \subseteq \Set{0}
 \cup \Set{\chi_E*\widetilde{\chi_E} = 0},
\eeq
since $\chi_E*\widetilde{\chi_E}$ is the Fourier transform of
$\Abs{\ft{\chi_E}}^2$
(where $\widetilde{g}(x) = \overline{g(-x)}$).
But
$$
\Set{\chi_E*\widetilde{\chi_E} = 0} = \Set{x:\ \Abs{E \cap (E+x)} = 0}.
$$
This and Lemma \ref{lm:shifts} imply that
$$
\supp{\ft{\delta_\Lambda}} \cap (-1,1) = \Set{0}.
$$

Let
$$
K_\delta(x) = \max\Set{0, 1-(1+\delta)\Abs{x}} = (1+\delta) \chi_{I_\delta}*\widetilde{\chi_{I_\delta}}(x),
$$
where $I_\delta = [0,{1\over 1+\delta}]$,
be a Fej\'er kernel (we will later take $\delta\to0$).
Then $\ft{K_\delta} = (1+\delta)\Abs{\ft{\chi_{I_\delta}}}^2$ is a
non-negative continuous function and,
after calculating $\ft{\chi_{I_\delta}}$, it follows that
$$
\ft{K_\delta}(0) = {1\over1+\delta}
$$
and
\beql{Kdzeros}
\Set{x: \ft{K_\delta}(x) = 0} = (1+\delta) ( \ZZ\setminus\Set{0} ).
\eeq

Next, we use the following result from \cite{K2}:
\begin{theorem}\label{th:at-zero}
Suppose that $\Lambda \in \RR^d$ is a multiset with density $\rho$,
$\delta_\Lambda = \sum_{\lambda\in\Lambda} \delta_\lambda$,
and that $\ft{\delta_\Lambda}$ is a measure in a neighborhood of
$0$.
Then $\ft{\delta_\Lambda}(\Set{0}) = \rho$.
\end{theorem}

\noindent
{\bf Remark.} The proof of Theorem \ref{th:at-zero}
shows that the assumption of $\ft{\delta_\Lambda}$ being a measure
in a neighborhood of zero is superfluous, if one knows a priori that
$\ft{\delta_\Lambda}$ is supported only at zero, in a neighborhood of zero.
Indeed, what is shown in that proof is that, as $t\to\infty$, the quantity
$\ft{\delta_\Lambda}(\phi(tx))$ remains bounded, for any $C_c^\infty$ test
function $\phi$. If $\ft{\delta_\Lambda}$ were not a measure
near $0$ but had support
only at $0$, locally, this quantity would grow like a polynomial in $t$ of degree
equal to the degree of the distribution at $0$.

Applying Theorem \ref{th:at-zero} and the remark following it we obtain that
$\ft{\delta_\Lambda}$ is equal to $\delta_0$ in a neighborhood of $0$, since $\Lambda$ has
density $1$.

Next, we claim that
$$
\sum_{\lambda\in\Lambda}\widehat{K_\delta}(x-\lambda)=1, \mbox{\ \ for a.e.\ $x$}.
$$
Indeed, take $\psi_\epsilon$ to be a smooth, positive-definite approximate
identity, supported in $(-\epsilon,\epsilon)$, and take $\epsilon = \epsilon(\delta)$
to be small enough so that $\supp{\psi_\epsilon*K_\delta} \subset (-1,1)$.
We have then
\begin{eqnarray*}
\sum_{\lambda\in\Lambda}\widehat{K_\delta}(x-\lambda)
 &=& \lim_{\epsilon\to 0} \sum_{\lambda\in\Lambda} \ft{\psi_\epsilon}(x-\lambda)
					\ft{K_\delta}(x-\lambda) \\
 &=& \lim_{\epsilon\to 0} \delta_\Lambda\left( (\ft{\psi_\epsilon}\ft{K_\delta})(x) \right) \\
 &=& \lim_{\epsilon\to 0} \ft{\delta_\Lambda}\left( (\psi_\epsilon*K_\delta)(x)\right)\\
 &=& \lim_{\epsilon\to 0} \delta_0 \left( (\psi_\epsilon*K_\delta)(x)\right)
     \mbox{\ \ \ (for $\epsilon$ small enough)}\\
 &=& \lim_{\epsilon\to 0} \psi_\epsilon*K_\delta(0) \\
 &=& K_\delta(0) \\
 &=& 1,
\end{eqnarray*}
which establishes the claim.
Applying this for $x\to 0$ and isolating the term $\lambda=0$ we get
$$
1 = {1\over 1+\delta} + \sum_{0\neq\lambda\in\Lambda} \ft{K_\delta}(-\lambda).
$$
Letting $\delta \to 0$ we obtain that $\ft{K_\delta}(-\lambda) \to 0$ for
each $\lambda\in\Lambda\setminus\Set{0}$, which implies that each such $\lambda$
is an integer, as
$\ZZ\setminus\Set{0}$ is the limiting set of the zeros of $\ft{K_\delta}$.

To get that $\Lambda = \ZZ$ notice that $\chi_{[0,1]}+\Lambda$
is a packing. By Theorem \ref{th:tiling-from-packing} again we get
that $\chi_{[0,1]}+\Lambda$ is in fact a tiling,
hence $\Lambda = \ZZ$.
\Qed


\section{Planar regions}


{\bf Proof of Theorem \ref{thm.2d}.}
We denote the coordinates in $\rr^2$ by $(x_1,x_2)$.
For $0\leq a\leq b\leq 1$ we will denote
\[
E_1(a,b)=(E\cap\{a\leq x_1\leq b, \ x_2\leq 0\})
\cup\{a\leq x_1\leq b, \ x_2\geq 0\},\]
\[
E_2(a,b)=(E\cap\{a\leq x_1\leq b, \ x_2\geq 0\})
\cup\{a\leq x_1\leq b, \ x_2\leq 0\},
\]
\[ F_1(a,b)=(E\cap\{a\leq x_2\leq b, \ x_1\leq 0\})
\cup\{a\leq x_2\leq b, \ x_1\geq 0\},
\]
\[ F_2(a,b)=(E\cap\{a\leq x_2\leq b, \ x_1\geq 0\})
\cup\{a\leq x_2\leq b, \ x_1\leq 0\}.
\]
We will also use $S_{a,b}$ to denote the vertical strip $[a,b]\times \RR$.
Let $v=(v_1,v_2)\in\RR^2$.
We will say that $E_2(a,b)$ {\em complements} $E_1(a',b')+v$ if
$E_1(a',b')+v$ is positioned above $E_2(a,b)$ so that (up to sets of measure
0) the two sets are disjoint and their union is $S_{a,b}$.
In particular, we must have $a'+v_1=a$ and $b'+v_1=b$. We will write
$\widetilde E_1(a,b)=S_{a,b}\setminus E_1(a,b)$, and similarly for $E_2$.
Finally, we write $A\sim B$ if the sets $A$ and $B$ are equal up to sets
of measure 0.

\begin{lemma}
Let $0<s''<s'<s<2s''$.
Suppose that $E_1(a,a+s)+v$, $E_1(a,a+s')+v'$, $E_1(a,a+s'')+v''$
complement $E_2(b-s,b)$, $E_2(b-s',b)$, $E_2(b-s'',b)$ respectively. Then
the points $v,v',v''$ are collinear. Moreover, the absolute value of the
slope of the line through $v,v''$ is bounded by $\eps(2s''-s)^{-1}$.
\label{2d.lemma}
\end{lemma}

Applying the lemma to the symmetric reflection of $E$ about the line
$x_2=1/2$, we find that the conclusions of the lemma also hold if we
assume that $E_2(a,a+s)+v$, $E_2(a,a+s')+v'$, $E_2(a,a+s'')+v''$
complement $E_1(b-s,b)$, $E_1(b-s',b)$, $E_1(b-s'',b)$ respectively.
Furthermore, we may interchange the $x_1$ and $x_2$
coordinates and obtain the analogue of the lemma with $E_1,E_2$ replaced
by $F_1,F_2$.

\medskip
{\bf Proof of Lemma \ref{2d.lemma}.}
Let $v=(v_1,v_2)$, $v'=(v'_1,v'_2)$, $v''=(v''_1,v''_2)$.
We first observe that if $v_1=v''_1$, it follows from the assumptions that
$v=v''$ and there is nothing to prove.  We may therefore assume
that $v_1\neq v''_1$.  We do, however, allow $v'=v$ or $v'=v''$.

\ignore{
We will write that $A\sim B$ if there are $c,c'>0$
such that (again up to sets of measure 0)
\[ A\cap \{x_2\leq c\}=B\cap\{x_2\leq c'\}.\]
If $A=E_1(a,b)$ and $B=E_1(a',b')+v$, $A\sim B$ means that the lower
parts of $E_1(a,b)$ and $E_1(a',b')$ coincide up to translation by $v$.
}

It follows from the assumptions that $E_2(b-s'',b)$ complements each of
$E_1(a,a+s'')+v''$, $E_1(a+s'-s'',a+s')+v'$, $E_1(a+s-s'',a+s)+v$.  Hence
\[
E_1(a+s'-s'',a+s')\sim E_1(a,a+s'')+(v''-v'),
\]
\[
E_1(a+s-s'',a+s)\sim E_1(a,a+s'')+(v''-v).
\]

Let $n$ be the unit vector perpendicular to $v-v''$ and such that
$n_2>0$.  For $t\in\rr$, let $P_t=\{x:\ x\cdot n\leq t\}$. We
define for $0\leq c\leq c'\leq 1$:
\[
\alpha_{c,c'}=\inf\{t\in\RR:\ |E_1(c,c')\cap P_t|>0\},
\]
\[
\beta_{c,c'}=\sup\{t\in\RR:\ |\widetilde E_1(c,c')\setminus P_t|>0\}.
\]
We will say that $x$ is a {\em low point} of $E_1(c,c')$  if
$x\in S_{c,c'}$, $x\cdot n=\alpha_{c,c'}$, and for any open disc
$D$ centered at $x$ we have
\begin{equation}
|D\cap E_1(c,c')|>0.
\label{e.low}\end{equation}
Similarly, we call $y$ a {\em high point} of $\widetilde E_1(c,c')$  if
$y\in S_{c,c'}$, $y\cdot n=\beta_{c,c'}$, and for any open disc
$D$ centered at $y$ we have
\begin{equation}
|D\cap \widetilde E_1(c,c')|>0.
\label{e.high}\end{equation}

It is easy to see that such points $x,y$ actually exist. Indeed,
by the definition of $\alpha_{c,c'}$ and an obvious covering
argument, for any $\alpha>\alpha_{c,c'}$
there are points $x'$ such that $x'\cdot n\leq \alpha$ and
that (\ref{e.low}) holds for any disc $D$ centered at $x'$.
Thus the set of such points $x'$ has at least one accumulation
point $x$ on the line $x\cdot n=\alpha_{c,c'}$.  It follows that
any such $x$ is a low point of $E_1(c,c')$.  The same argument
works for $y$.

The low and high points need not be unique;
however, all low points $x$ of $E_1(c,c')$ lie on the same
line $x\cdot n=\alpha_{c,c'}$ parallel to the vector $v-v''$,
and similarly for high points. Furthermore, the
low and high points of $E_1(c,c')$ do not change if $E_1(c,c')$
is modified by a set of measure 0.

Let now $A= E_1(a,a+s'')$, and let $x$ be a low point of
$A$.  Since $s<2s''$, we have
\[
B:= E_1(a,a+s)= E_1(a,a+s'')\cup E_1(a+s-s'',a+s)
\sim A\cup (A+v''-v),
\]
hence $x$ is also a low point of $B$ with respect to $v-v''$.
Now note that
\[
E_1(a+s'-s'',a+s')\sim A+(v''-v')
\]
intersects any open neighbourhood of $x+(v''-v')$ in positive
measure.  But on the other hand, $E_1(a+s'-s'',a+s')
\subset B$.  By the extremality of $x$ in $B$, $x+(v''-v')$ lies
on or above the line segment joining $x$ and $x+(v''-v)$, hence
$v''-v'$ lies on or above the line segment joining $0$ and $v''-v$.

Repeating the argument in the last paragraph with $x$ replaced by
a high point $y$ of $\widetilde E_1(a,a+s'')$, we obtain that $v''-v'$ lies on or below
the line segment joining $0$ and $v''-v$. Hence $v,v',v''$ are collinear.

Finally, we estimate the slope of the line through $v,v''$.
We have to prove that
\begin{equation}
\frac{2s''-s}{s-s''}\,|v''_2-v_2|\leq\eps
\label{slope}\end{equation}
(recall that $v''_1-v_1=s-s''$).
Define $x$ as above, and let $k\in\zz$.  Iterating translations
by $v-v''$ (in both directions), we find that $x+k(v-v'')$ is
a low point of $B$ as long as it belongs to $B$, i.e. as long as
\[ a\leq x_1+k(s-s'')\leq a+s.\]
The number of such $k$'s is at least $\frac{s}{s-s''}-1$.
On the other hand, all low points of $B$ lie in the rectangle
$a\leq x_1\leq a+s, -\eps\leq x_2\leq 0$. Hence
\[ (\frac{s}{s-s''}-2)|v''_2-v_2|\leq\eps,\]
which is (\ref{slope}).
\Qed

\begin{figure}[htbp] \centering
\ \psfig{figure=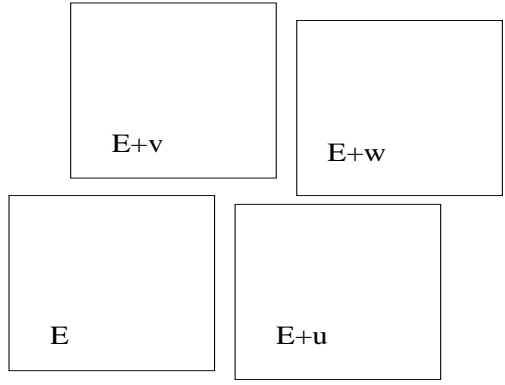,height=2in,width=2.6in}
\caption{A ``corner" and a fourth near-square.}
        \label{fig:corner}
        \end{figure}

\bigskip
We return to the proof of Theorem \ref{thm.2d}.  Since $E$ is almost
a square, we know roughly how the translates of $E$ can fit together.
Locally, any tiling by $E$ is essentially a tiling by a ``solid" $1\times 1$
square with ``margins" of width between $0$ and $2\eps$ (see Fig.
\ref{fig:corner}).

We first locate a ``corner".  Namely, we may assume
that the tiling contains $E$ and its translates $E+u$, $E+v$, where
\begin{equation}\label{e.u}
1\leq u_1\leq 1+2\eps,\ -2\eps\leq u_2\leq 2\eps,
\end{equation}
\begin{equation}\label{e.v}
0\leq v_1\leq \half+\eps,\ 1\leq v_2\leq 1+2\eps.
\end{equation}
This can always be achieved by translating the tiled plane and
taking symmetric reflections of it if necessary.

Let $E+w$ be the translate of $E$ which fits into this corner:
\begin{equation}\label{e.w}
v_1+1\leq w_1\leq v_1+1+2\eps,\ u_2+1\leq w_2\leq u_2+1+2\eps.
\end{equation}
We will prove that $w=u+v$ (without the $\eps$-errors).

\relax From (\ref{e.w}), (\ref{e.u}), (\ref{e.v}) we have
\[
1\leq w_1\leq \frac{3}{2}+3\eps,\ -4\eps\leq w_2-v_2\leq 4\eps.\]
Hence $w$ satisfies both of the following.

\medskip\noindent
(A) $E_2(0,1-(w_1-u_1))$ complements $E_1(w_1-u_1,1)+(w-u)$, and
\[ 1-(w_1-u_1)=1-w_1+u_1\geq 1+1-(\frac{3}{2}+3\eps)=\half -3\eps,\]
\[ |(w_1-u_1)-v_1|=|(w_1-v_1)-u_1|\leq 2\eps.\]

\medskip \noindent
(B) $-4\eps\leq w_2-v_2\leq 4\eps$, and $F_2(r,t)$ complements
$F_1(r',t')+(w-v)$, where
\[ r=\max(0,w_2-v_2), \ r'=\max(0,v_2-w_2),\]
\[ t=1-\max(0,v_2-w_2), \ t'=1-\max(0,w_2-v_2).\]

\medskip
If $w=u+v$, we have $w-u=v$, $w-v=u$, hence by considering the
``corner" $E,E+u,E+v$ we see that both (A) and (B) hold.
Assuming that $\eps$ is small enough, we shall prove that:

\medskip
1. All points $w$ satisfying (A) lie on a fixed straight line $l_1$
making an angle less than $\pi/4$ with the $x_1$ axis.

\medskip
2. All points $w$ satisfying (B) lie on a fixed straight line $l_2$
making an angle at most $\pi/4$ with the $x_2$ axis.

\medskip
It follows that there can be at most one $w$ which satisfies both
(A) and (B), since $l_1$ and $l_2$ intersect only at one point.
Consequently, if $E+w$ is the
translate of $E$ chosen as above, we must have $w=u+v$.
Now it is easy to see that $E+\Lambda$ is a tiling, where $\Lambda$
is the lattice $\{ku+mv:\ k,m\in\ZZ\}$.

We first prove 1.  Suppose that $w,w', w'',\dots$ (not necessarily
all distinct) satisfy (A). By the assumptions in (A), we may apply Lemma
\ref{2d.lemma} with $E_1$ and $E_2$ interchanged and with $a=0$,
$b=1$, $s=1-(w_1-u_1),s'=1-(w'_1-u_1),\dots\geq \half-3\eps$.
\relax From the second inequality in (A) and the triangle inequality
we also have $|s-s''|\leq 4\eps$.  We find that all $w$ satisfying (A)
lie on a line $l_1$ with slope bounded by
\[
\frac{\eps}{|2s''-s|}\leq \frac{\eps}{s''-|s''-s|}
\leq \frac{\eps}{1/2-7\eps},\]
which is less than 1 if $\eps< 1/16$.

To prove 2., we let $w,w',w''$ be three (not necessarily distinct)
points satisfying (B) and such that $w_2\leq w'_2\leq w''_2$. We then
apply the obvious analogue of Lemma \ref{2d.lemma} with $E_1,E_2$
replaced by $F_1,F_2$ and with $a=\max(v_2-w_2,0)\leq 4\eps$,
$b=1-\max(v_2-w_2)\geq 1-4\eps$.  From the estimates in (B) we have
$1-16\eps\leq s,s',s''\leq 1$, hence
$|2s''-s|\geq 2-32\eps-1=1-32\eps$.
We conclude that all $w$ satisfying (B) lie on a line $l_2$ such
that the inverse of the absolute value of its slope is bounded by
$\frac{\eps}{1-32\eps}$. This is at most 1 if $\eps \leq 1/33$.
\Qed

\bigskip
{\bf Proof of Corollary \ref{cor.2d}.}  Let $Q=[0,1]\times[0,1]$.
By rescaling, it suffices to prove that for any $\eps>0$ there is a
$\delta>0$ such that if $E\subset Q$, $E$ tiles $\RR^2$ by translations,
and $|E|\geq 1-\delta$, then $E$ contains the square $$Q_\eps= [\eps,1-\eps]
\times [\eps,1-\eps]$$ (up to sets of measure 0). The result then
follows from Theorem \ref{thm.2d}.

Let $E$ be as above, and suppose that $Q_\eps\setminus E$
has positive measure.
Since $E$ tiles $\RR^2$, there is
a $v\in \RR^2$ such that $|E\cap (E+v)|=0$ and $|Q_\eps\cap
(E+v)|>0$.  We then have
\[
|E\cup (E+v)|=|E|+|E+v|\geq 2-2\delta,
\]
but also
\[
|E\cup (E+v)|\leq |Q\cup (Q+v)|\leq 2-\eps^2,
\]
since $E\subset Q$, $E+v\subset Q+v$, and $Q_\eps\cap (Q+v)\neq
\emptyset$ so that $|Q\cap (Q+v)|\geq\eps^2$.  This is 
a contradiction if $\delta$ is small enough.
\Qed


\section{A counterexample in higher dimensions}
\label{sec-3d}


In this section we prove Theorem \ref{3d-thm}.  It suffices to construct
$E$ for $n=3$, since then $E\times [0,1]^{n-3}$ is a subset of $\RR^n$
with the required properties.

Let $(x_1,x_2,x_3)$ denote the Cartesian coordinates in $\rr^3$.  
It will be convenient to rescale $E$ so that 
$[\eps,1]^3\subset E\subset [0,1+\eps]^3$.

\begin{figure}[htbp] \centering
\ \psfig{figure=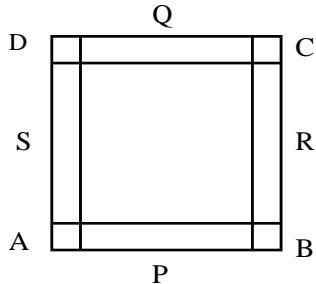,height=1.5in,width=1.6in}
\caption{The construction of $E$.}
        \label{fig:divide}
        \end{figure}

We construct $E$ as follows.  We let $E$ be bounded from below and above by the
planes $x_3=0$ and $x_3=1$ respectively.  The planes $x_1=\eps, x_1=1,x_2=\eps,x_2=1$
divide the cube $[0,1+\eps]^3$ into 9 parts (Figure \ref{fig:divide}).  The
middle part is entirely contained in $E$.  We label by $A,B,C,D,P,Q,R,S$ the
remaining 8 segments as shown in Figure \ref{fig:divide}.  We then let
\[
E\cap P= P\cap\Big\{0\leq x_3\leq \frac{1}{8}\hbox{ or }
\frac{1}{2}\leq x_3\leq \frac{5}{8}\Big\},
\]
\[
E\cap R= R\cap\Big\{0\leq x_3\leq \frac{1}{8}\hbox{ or }
\frac{1}{2}\leq x_3\leq \frac{5}{8}\Big\},
\]
\[
E\cap Q= Q\cap\Big\{0\leq x_3\leq \frac{1}{4}\hbox{ or }\frac{3}{8}\leq x_3\leq \frac{3}{4}
\hbox{ or }\frac{7}{8}\leq x_3\leq 1\Big\},
\]
\[
E\cap S= S\cap\Big\{0\leq x_3\leq \frac{1}{4}\hbox{ or }\frac{3}{8}\leq x_3\leq \frac{3}{4}
\hbox{ or }\frac{7}{8}\leq x_3\leq 1\Big\},
\]
and
\[
E\cap A= A\cap\Big\{0\leq x_3\leq \frac{1}{16}\Big\},\ 
E\cap C= A\cap\Big\{\frac{1}{2}\leq x_3\leq \frac{9}{16}\Big\},\ 
\]
\[
E\cap B= B\cap\Big\{\frac{5}{16}\leq x_3\leq \frac{3}{4}\Big\},\ 
E\cap D= D\cap\Big\{0\leq x_3\leq\frac{1}{4}\hbox{ or }\frac{13}{16}\leq x_3\leq 1\Big\}.
\]
We also denote $K=\bigcup_{j\in\zz}(E+(0,0,j))$.

Let $E+T$ be a tiling of $\RR^3$, and assume that $0\in T$.  Suppose that 
$E+v$ and $E+w$ are neighbours 
in this tiling so that the vertical sides of $(E\cap P)+v$ and $(E\cap Q)+w$ 
meet in a set of non-zero two-dimensional measure.  Then we must have
$v-w=(0,1,(v-w)_3)$, where $(v-w)_3\in\{\pm\frac{1}{4},\pm\frac{3}{4}\}$.
A similar statement holds with $P,Q$ replaced by $R,S$ and with the $x_1,x_2$
coordinates interchanged.
We deduce that the tiling consists of copies of $E$ stacked into identical vertical
``columns" $K_{ij}=K+(i,j,t_{ij})$, arranged in a rectangular grid in the $x_1x_2$
plane and shifted vertically so that $t_{i+1,j}-t_{ij}$ and 
$t_{i,j+1}-t_{ij}$ are always $\pm \frac{1}{4}$.  We will use matrices 
$(t_{ij})$ to encode such a tiling or portions thereof.

It is easy to see that $(t_{ij})$, where $t_{ij}=0$ if $i+j$ is even
and $\frac{1}{4}$ if $i+j$ is odd, is indeed a tiling.  
It remains to show that $E$ does not admit a lattice tiling.
Indeed, the four possible choices of the generating vectors  in
any lattice $(t_{ij})$ with $t_{ij}=\pm\frac{1}{4}$ produce the configurations
\[
\pmatrix{0&t\cr t&2t\cr},\ 
\pmatrix{2t&t\cr t&0\cr},\ 
\pmatrix{0&t\cr -t&0\cr},\ 
\pmatrix{0&-t\cr t&0\cr}.
\]
But it is easy to see that the corners $A,B,C,D$ do not match if so translated.



{\sc Department of Mathematics, University of Crete, Knossos Ave.,
714 09 Iraklio, Greece.}
E-mail: {\tt mk@fourier.math.uoc.gr}

{\sc Department of Mathematics, University of British Columbia,
Vancouver, B.C. V6T 1Z2, Canada.}
E-mail: {\tt ilaba@math.ubc.ca}


\begin{thebibliography}{99}

\bibitem{Fug} B. Fuglede: {\it Commuting self-adjoint partial differential
operators and a group-theoretic problem}, J. Funct. Anal. 16 (1974),
101--121

\bibitem{GN} D. Girault-Beauquier, M. Nivat: {\it Tiling the plane
with one tile}, in: {\it Topology and Category Theory in Computer
Science}, G.M. Reed, A.W. Roscoe, R.F. Wachter (eds.), Oxford
Univ. Press 1989, 291--333.

\bibitem{GS} B. Gr\"unbaum, G.C. Shepard: {\it Tilings and patterns},
New York: Freeman 1987.

\bibitem{IKP} A. Iosevich, N. H. Katz, S. Pedersen: {\it Fourier bases and a
distance problem of Erd\"os}, Math. Res. Letters 6 (1999), 251--255.

\bibitem{IKT1} A. Iosevich, N. H. Katz, T. Tao: {\it Convex bodies with a
point of curvature do not have Fourier bases}, Amer. J. Math. 123 (2001),
115--120.

\bibitem{IKT2} A. Iosevich, N.H. Katz, T.Tao: {\it Fuglede conjecture
holds for convex planar domains}, preprint, 2001.

\bibitem{JP1} P. Jorgensen, S. Pedersen: {\it Spectral pairs in
Cartesian coordinates},  J. Fourier Anal. Appl. 5 (1999), 285--302.

\bibitem{Ke} R. Kenyon: {\it Rigidity of planar tilings}, Invent. Math.
107 (1992), 637--651.

\bibitem{K3} M. Kolountzakis: {\it Non-symmetric convex domains have no
basis of exponentials}, Illinois J. Math. 44 (2000), 542--550.

\bibitem{K1} M.N. Kolountzakis: {\it Packing, Tiling, Orthogonality and Completeness,}
Bull.\ L.M.S. {\bf 32} (2000), 5, 589--599.

\bibitem{K2} M.N. Kolountzakis: {\it On the structure of multiple translational tilings by
polygonal regions}, Discrete Comput.\ Geom.\ {\bf 23} (2000), 4, 537--553.

\bibitem{KL} M.N. Kolountzakis, J.C. Lagarias: {\it Structure of tilings of the
line by a function}, Duke Math. J. 82 (1996), 653--678.

\bibitem{KP} M.N. Kolountzakis, M. Papadimitrakis: {\it A class of non-convex
polytopes that admit no orthonormal basis of exponentials}, preprint, 2001.

\bibitem{L1} I. {\L}aba: {\it Fuglede's conjecture for a union of two intervals},
Proc. AMS 121 (2001), 2965--2972. 

\bibitem{L2} I. {\L}aba: {\it The spectral set conjecture and multiplicative
properties of roots of polynomials}, J. London Math. Soc., to appear.

\bibitem{LW1} J. Lagarias, Y. Wang: {\it Tiling the line with translates
of one tile}, Inv. Math. 124 (1996), 341--365.

\bibitem{LW2} J. Lagarias, Y. Wang: {\it Spectral sets and factorization of
finite abelian groups}, J. Funct. Anal. 73 (1997), 122--134.

\bibitem{M} P. McMullen: {\it Convex bodies which tile the space by translation},
Mathematika 27 (1980), 113--121.

\bibitem{V} B.A. Venkov: {\it On a class of Euclidean polyhedra}, Vestnik
Leningrad Univ. Ser. Mat. Fiz. Him. 9 (1954), 11--31.
\end{thebibliography}
\end{document}